\documentclass[12pt]{amsart}%
\usepackage{amsmath}
\usepackage{graphicx}
\usepackage{amsfonts}
\usepackage{amssymb}%
\newtheorem{theorem}{Theorem}
\theoremstyle{plain}

\newtheorem{lemma}{Lemma}

\newtheorem{remark}{Remark}

\numberwithin{equation}{section}
\begin{document}
\title[Nonlinear pattern formation]{Pattern formation (II): The Turing Instability}
\author{Yan Guo}
\address{Division of Applied Mathematics, Brown University, Providence, RI 02912, USA}
\email{guoy@dam.brown.edu}
\author{Hyung Ju Hwang}
\address{School of Mathematics, Trinitiy College Dublin\\
Dublin 2, Ireland}
\email{hjhwang@maths.tcd.ie}
\date{}

\begin{abstract}
We consider the classical Turing instability in a reaction-diffusion system as
the secend part of our study on pattern formation. We prove that nonlinear
dynamics of a general perturbation of the Turing instability is determined by
the finite number of linear growing modes over a time scale of $\ln\frac
{1}{\delta},$ where $\delta$ is the strength of the initial perturbation.

\end{abstract}
\maketitle

\section{Growing modes in a reaction-diffusion system}

In this section we summarize the classical linear Turing instability criterion
for a reaction-diffusion system. Consider a reaction-diffusion system of
$2$-species as%
\begin{align}
\frac{\partial U}{\partial t}  &  =\nabla\cdot\left(  D_{1}\left(
U\mathbf{,}V\right)  \nabla U\right)  +f\left(  U\mathbf{,}V\right)
,\label{RD}\\
\frac{\partial V}{\partial t}  &  =\nabla\cdot\left(  D_{2}\left(
U\mathbf{,}V\right)  \nabla V\right)  +g\left(  U\mathbf{,}V\right)
,\nonumber
\end{align}
where $U\left(  \mathbf{x,}t\right)  \mathbf{,}V\left(  \mathbf{x,}t\right)  $
are concentration for species, $D_{1},D_{2}\,$\ diffusion coefficients, $f,g$
reaction terms.$\ $

In this paper we consider a $d$-dimensional box $\mathbb{T}^{d}=\left(
0,\pi\right)  ^{d},$ $d=1,2,3,$ with Neumann boundary conditions for $U$ and
$V$, i.e.,$\ $%
\begin{equation}
\frac{\partial U}{\partial x_{i}}=\frac{\partial V}{\partial x_{i}}=0\text{
\ at }x_{i}=0,\pi,\ \ \text{for }1\leq i\leq d. \label{Neumann}%
\end{equation}
Homogeneous steady state $U=\bar{U}\mathbf{,}V\mathbf{=}\bar{V}$ forms a
steady state provided%
\begin{equation}
0=f\left(  \bar{U}\mathbf{,}\bar{V}\right)  =g\left(  \bar{U}\mathbf{,}\bar
{V}\right)  . \label{steady-state}%
\end{equation}
In this article, we study the nonlinear evolution of a perturbation
\[
u(x,t)=U(x,t)-\bar{U},\;v(x,t)=V(x,t)-\bar{V}%
\]
around $[\bar{U},\bar{V}]$, which satisfies the equivalent reaction-diffusion
system:%
\begin{align}
\frac{\partial u}{\partial t}  &  =\nabla\cdot\left(  D_{1}\left(  u+\bar
{U}\mathbf{,}v+\bar{V}\right)  \nabla u\right)  +f\left(  u+\bar{U}%
\mathbf{,}v+\bar{V}\right)  ,\label{nonlinear1}\\
\frac{\partial v}{\partial t}  &  =\nabla\cdot\left(  D_{2}\left(  u+\bar
{U}\mathbf{,}v+\bar{V}\right)  \nabla v\right)  +g\left(  u+\bar{U}%
\mathbf{,}v+\bar{V}\right)  . \label{nonlinear2}%
\end{align}
The corresponding linearized system then takes the form%
\begin{align}
u_{t}  &  =\bar{D}_{1}\nabla^{2}u+\bar{f}_{u}u+\bar{f}_{v}v,\label{linear1}\\
v_{t}  &  =\bar{D}_{2}\nabla^{2}v+\bar{g}_{u}u+\bar{g}_{v}v, \label{linear2}%
\end{align}
where $\bar{D}_{1}=D_{1}\left(  \bar{U},\bar{V}\right)  ,\ \bar{D}_{2}%
=D_{2}\left(  \bar{U},\bar{V}\right)  ,\ \bar{f}_{u}=\frac{\partial
f}{\partial u}\left(  \bar{U},\bar{V}\right)  ,\bar{f}_{v}=\frac{\partial
f}{\partial v}\left(  \bar{U},\bar{V}\right)  ,\bar{g}_{u}=\frac{\partial
g}{\partial u}\left(  \bar{U},\bar{V}\right)  ,\bar{g}_{v}=\frac{\partial
g}{\partial v}\left(  \bar{U},\bar{V}\right)  $.

We use $[\cdot,\cdot]$ to denote a column vector, and let
\[
\mathbf{w}(x,t)\equiv\lbrack u(x,t),v(x,t)],\ \mathbf{\bar{W}=[}\bar{U}%
,\bar{V}].
\]
Then the original nonlinear system (\ref{nonlinear1}) and (\ref{nonlinear2})
can be written in a matrix form:%
\begin{align}
\frac{\partial\mathbf{w}}{\partial t}=  &  \nabla\cdot\left(  D\nabla
\mathbf{w}\right)  \mathbf{+F}\label{matrix}\\
=  &  (\bar{D}\nabla^{2}\mathbf{w+}A\mathbf{w)}\mathbf{+(\{}\nabla\cdot\left(
D\nabla\mathbf{w}\right)  -\bar{D}\nabla^{2}\mathbf{w\}+F-}A\mathbf{w)}%
\nonumber\\
\equiv &  \mathcal{L}\left(  \mathbf{w}\right)  +\mathcal{N}\left(
\mathbf{w}\right)  .\nonumber
\end{align}
where%
\begin{align*}
D  &  =\left(
\begin{array}
[c]{cc}%
D_{1}\left(  \mathbf{w+\bar{W}}\right)  & 0\\
0 & D_{2}\left(  \mathbf{w+\bar{W}}\right)
\end{array}
\right)  ,\ \bar{D}=\left(
\begin{array}
[c]{cc}%
\bar{D}_{1} & 0\\
0 & \bar{D}_{2}%
\end{array}
\right)  ,\\
\mathbf{F}  &  \mathbf{=}\left(
\begin{array}
[c]{c}%
f\left(  \mathbf{w+\bar{W}}\right) \\
g\left(  \mathbf{w+\bar{W}}\right)
\end{array}
\right)  ,\ A=\left(
\begin{array}
[c]{cc}%
\bar{f}_{u} & \bar{f}_{v}\\
\bar{g}_{u} & \bar{g}_{v}%
\end{array}
\right)  .
\end{align*}
Let $\mathbf{q=}\left(  q_{1},..,q_{d}\right)  \in\Omega=\left(
\mathbb{N}\cup\left\{  0\right\}  \right)  ^{d}$ and let%
\[
e_{\mathbf{q}}(x)\equiv\prod_{i=1}^{d}\cos\left(  q_{i}x_{i}\right)  ,
\]
where $\mathbf{q\in}\Omega$. Then $\left\{  e_{\mathbf{q}}(x)\right\}
_{\mathbf{q\in}\Omega}$ forms a basis of the space of functions in
$\mathbb{T}^{d}$ that satisfy Neumann boundary condition (\ref{Neumann}).

We look for a normal mode to the linear reaction-diffusion system
(\ref{linear1}) and (\ref{linear2}) of the following form:%
\begin{equation}
\mathbf{w}\left(  x,t\right)  =\mathbf{r}_{\mathbf{q}}\exp\left(
\lambda_{\mathbf{q}}t\right)  e_{\mathbf{q}}(x), \label{normal-mode}%
\end{equation}
where $\mathbf{r}_{\mathbf{q}}$ is a vector depending on $\mathbf{q.}$ We
substitute (\ref{normal-mode}) into (\ref{linear1})-(\ref{linear2}) to get%
\[
\lambda_{\mathbf{q}}\mathbf{r}_{\mathbf{q}}=\left(
\begin{array}
[c]{cc}%
\bar{f}_{u}-\bar{D}_{1}q^{2} & \bar{f}_{v}\\
\bar{g}_{u} & \bar{g}_{v}-\bar{D}_{2}q^{2}%
\end{array}
\right)  \mathbf{r}_{\mathbf{q}},
\]
where $q^{2}=\sum_{i=1}^{d}q_{i}^{2}$. A nontrivial normal mode can be
obtained by setting%
\[
\det\left(
\begin{array}
[c]{cc}%
\lambda_{\mathbf{q}}-\bar{f}_{u}+\bar{D}_{1}q^{2} & -\bar{f}_{v}\\
-\bar{g}_{u} & \lambda_{\mathbf{q}}-\bar{g}_{v}+\bar{D}_{2}q^{2}%
\end{array}
\right)  =0.
\]
This leads to the following dispersion formula for $\lambda_{\mathbf{q}}$:%
\begin{equation}
\lambda_{\mathbf{q}}^{2}+\{-\bar{f}_{u}+\bar{D}_{1}q^{2}-\bar{g}_{v}+\bar
{D}_{2}q^{2}\}\lambda_{\mathbf{q}}+\{\left(  \bar{f}_{u}-\bar{D}_{1}%
q^{2}\right)  \left(  \bar{g}_{v}-\bar{D}_{2}q^{2}\right)  -\bar{f}_{v}\bar
{g}_{u}\}=0. \label{dispersion}%
\end{equation}
We assume first that without diffusion, the $\lambda_{\mathbf{q}}$ has
negative real part (stable):
\begin{equation}
\text{tr }A=\bar{f}_{u}+\bar{g}_{v}<0,\ \det A=\bar{f}_{u}\bar{g}_{v}-\bar
{f}_{v}\bar{g}_{u}>0, \label{tr-det}%
\end{equation}
On the other hand, in the presence of diffusion, we assume the following
diffusion-driven (linear) instability criterion by requiring there exists a
$q$ such that%
\begin{equation}
\left(  \bar{f}_{u}-\bar{D}_{1}q^{2}\right)  \left(  \bar{g}_{v}-\bar{D}%
_{2}q^{2}\right)  -\bar{f}_{v}\bar{g}_{u}<0, \label{instability-criteria}%
\end{equation}
which ensures that (\ref{dispersion}) has at least one positive root
$\lambda_{\mathbf{q}}.$

\begin{remark}
To satisfy (\ref{tr-det}) and (\ref{instability-criteria}), the discriminant
for the quadratic equation for $q^{2}$ in (\ref{instability-criteria}) must be
positive:
\begin{equation}
\left(  \bar{f}_{u}\bar{D}_{2}+\bar{g}_{v}\bar{D}_{1}\right)  >2\sqrt{\bar
{D}_{1}\bar{D}_{2}}\det A>0, \label{range}%
\end{equation}
which means the range of inhibition $\sqrt{\bar{D}_{2}/\left\vert \bar{g}%
_{v}\right\vert }$ is larger than the range of activation $\sqrt{\bar{D}%
_{1}/\left\vert \bar{f}_{u}\right\vert }.$ From (\ref{tr-det}) and
(\ref{range}), it follows that
\begin{equation}
\bar{f}_{u}\bar{g}_{v}<0,\ \ \text{and \ \ }\bar{f}_{v}\bar{g}_{u}<0,
\label{sign}%
\end{equation}
and we have only two cases for $A:$%
\[
A=\left(
\begin{array}
[c]{cc}%
+ & -\\
+ & -
\end{array}
\right)  \text{ \ or \ }A=\left(
\begin{array}
[c]{cc}%
+ & +\\
- & -
\end{array}
\right)  ,
\]
where formal case is called activator-inhibitor (or predator-prey) and the
latter positive feedback. It also follows from (\ref{tr-det}) that
\[
\bar{D}_{1}\neq\bar{D}_{2}.
\]

\end{remark}

For given $\mathbf{q}\in\Omega$, we denote the corresponding eigenvalues by
$\lambda_{\pm}(\mathbf{q})$ and eigenvectors by $\mathbf{r}_{\pm}(\mathbf{q}%
)$. We split into the three cases for the linear analysis:

(1) \ Generic case where we have two independent real eigenvectors and we
denote%
\[
\Omega_{\text{generic}}\equiv\{\mathbf{q}\in\Omega\text{ such that }%
\mathbf{r}_{+}(\mathbf{q})\neq\mathbf{r}_{-}(\mathbf{q})\}.
\]

By an elementary computation of the discriminant of (\ref{dispersion}), we
have, except for only \textit{finitely} many $q$,
\[
\left(  \bar{D}_{1}-\bar{D}_{2}\right)  q^{4}-\text{tr }A\ \left(  \bar{D}%
_{1}+\bar{D}_{2}\right)  q^{2}+4\left(  \bar{f}_{u}\bar{D}_{2}+\bar{g}_{v}%
\bar{D}_{1}\right)  q^{2}+\left(  \text{tr }A\right)  ^{2}-4\det A>0,
\]
since $\bar{D}_{1}-\bar{D}_{2}\neq0.$ Therefore, there are two distinct real
roots such that%
\[
\lambda_{-}(\mathbf{q})<\lambda_{+}(\mathbf{q})
\]
for large $q.$ Since $\bar{f}_{v}\neq0$ in (\ref{sign}), the corresponding
(linearly independent) eigenvectors $\mathbf{r}_{-}(\mathbf{q})$ and
$\mathbf{r}_{+}(\mathbf{q})$ are given by
\begin{equation}
\mathbf{r}_{\pm}(\mathbf{q})=\left[  1,\frac{\lambda_{\pm}(\mathbf{q})-\bar
{f}_{u}+\bar{D}_{1}q^{2}}{\bar{f}_{v}}\right]  .
\end{equation}
It is easy to see from (\ref{instability-criteria}) that there exist only
finitely many $\mathbf{q}$ such that $\lambda_{+}(\mathbf{q})>0.$ We therefore
can denote the largest eigenvalue by $\lambda_{\text{max}}>0$ and define
\[
\Omega_{\text{max}}\equiv\{\mathbf{q}\in\Omega\text{ such that }\lambda
_{+}(\mathbf{q})=\lambda_{\text{max }}\}.
\]
We also denote $\nu>0$ to be the gap between the $\lambda_{\text{max}}$ and
the rest. Moreover, there is one $q^{2}$ (possibly two) having $\lambda
_{\mathbf{q}}^{+}\left(  q^{2}\right)  =\lambda_{\text{max}}$ when we regard
$\lambda_{\mathbf{q}}^{+}$ as a function of $q^{2}$.

(2) Defective case where we have the repeated real eigenvalues and eigenvectors:

Note that there may be possibly one $q^{2}$ (so finitely many $\mathbf{q}$)
such that from (\ref{tr-det})
\begin{equation}
\lambda_{+}(\mathbf{q})=\lambda_{-}(\mathbf{q})\equiv\lambda(\mathbf{q}%
)=\{\bar{f}_{u}+\bar{g}_{v}-\left(  \bar{D}_{1}+\bar{D}_{2}\right)
q^{2}\}/2<0 \label{<0}%
\end{equation}
and$\ \mathbf{r}_{+}(\mathbf{q})=\mathbf{r}_{-}(\mathbf{q})\equiv
\mathbf{r}(\mathbf{q})$ and we denote%
\[
\Omega_{\text{defective}}\equiv\{\mathbf{q}\in\Omega\text{ such that
}\mathbf{r}_{+}(\mathbf{q})=\mathbf{r}_{-}(\mathbf{q})\}.
\]
In this case we find another independent vector
\[
\mathbf{r}^{\prime}(\mathbf{q})=[0,\frac{1}{\bar{f}_{v}}]
\]
satisfying $\left(  A-\lambda(\mathbf{q})I\right)  \mathbf{r}^{\prime
}(\mathbf{q})$ $=\mathbf{r}(\mathbf{q}).$

(3) Complex case where we have complex eigenvalues for $q$ and we denote it by
$\Omega_{\text{complex}}\equiv\Omega-\left(  \Omega_{\text{generic}}\cup
\Omega_{\text{defective}}\right)  $. For $\mathbf{q}\in\Omega_{\text{complex}%
},$ we denote $\lambda_{+}(\mathbf{q})\equiv\operatorname{Re}\lambda
(\mathbf{q})+i\operatorname{Im}\lambda(\mathbf{q})$ and $\mathbf{r}%
_{+}(\mathbf{q})\equiv\operatorname{Re}\mathbf{r}(\mathbf{q}%
)+i\operatorname{Im}\mathbf{r}(\mathbf{q})$. Then we have $\lambda
_{-}(\mathbf{q})\equiv\operatorname{Re}\lambda(\mathbf{q})-i\operatorname{Im}%
\lambda(\mathbf{q})$ and $\mathbf{r}_{-}(\mathbf{q})\equiv\operatorname{Re}%
\mathbf{r}(\mathbf{q})-i\operatorname{Im}\mathbf{r}(\mathbf{q})$. Notice that
$\operatorname{Re}\lambda(\mathbf{q})<0$ as in (\ref{<0}), and
$\operatorname{Re}\mathbf{r}(\mathbf{q})$ and $\operatorname{Im}%
\mathbf{r}(\mathbf{q})$ are linearly independent vectors.

Given any initial perturbation $\mathbf{w}\left(  \mathbf{x},0\right)  $, we
can expand it as%
\begin{align*}
\mathbf{w}\left(  \mathbf{x},0\right)   &  =\sum_{_{\mathbf{q}}\in\Omega
}\mathbf{w}_{\mathbf{q}}e_{\mathbf{q}}(x)=\sum_{_{\mathbf{q}}\in
\Omega_{\text{generic}}}\{w_{\mathbf{q}}^{-}\mathbf{r}_{-}(\mathbf{q}%
)+w_{\mathbf{q}}^{+}\mathbf{r}_{+}(\mathbf{q})\}e_{\mathbf{q}}(x)\\
&  +\sum_{_{\mathbf{q}}\in\Omega_{\text{defective}}}\mathbf{\{}w_{\mathbf{q}%
}\mathbf{r}(\mathbf{q})+w_{\mathbf{q}}^{\prime}\mathbf{r}^{\prime}%
(\mathbf{q})\}e_{\mathbf{q}}(x)\\
&  +\sum_{_{\mathbf{q}}\in\Omega_{\text{complex}}}\{w_{\mathbf{q}%
}^{\operatorname{Re}}\operatorname{Re}\mathbf{r}(\mathbf{q})+w_{\mathbf{q}%
}^{\operatorname{Im}}\operatorname{Im}\mathbf{r}(\mathbf{q})\}e_{\mathbf{q}%
}(x),
\end{align*}
so that
\begin{align}
\mathbf{w}_{\mathbf{q}}  &  =w_{\mathbf{q}}^{-}\mathbf{r}_{-}(\mathbf{q}%
)+w_{\mathbf{q}}^{+}\mathbf{r}_{+}(\mathbf{q})\text{ \ for \ \ }\mathbf{q}%
\in\Omega_{\text{generic}},\label{inde}\\
\mathbf{w}_{\mathbf{q}}  &  =w_{\mathbf{q}}\mathbf{r}(\mathbf{q}%
)+w_{\mathbf{q}}^{\prime}\mathbf{r}^{\prime}(\mathbf{q})\ \ \ \ \ \text{for
\ \ }\mathbf{q}\in\Omega_{\text{defective}},\nonumber\\
\mathbf{w}_{\mathbf{q}}  &  =w_{\mathbf{q}}^{\operatorname{Re}}%
\operatorname{Re}\mathbf{r}(\mathbf{q})+w_{\mathbf{q}}^{\operatorname{Im}%
}\operatorname{Im}\mathbf{r}(\mathbf{q})\text{ for \ \ }\mathbf{q}\in
\Omega_{\text{complex}}.\nonumber
\end{align}
The unique solution $\mathbf{w}\left(  x,t\right)  =[u\left(  x,t\right)
,v\left(  x,t\right)  ]$ to (\ref{linear1})-(\ref{linear2}) is given by
\begin{align}
\mathbf{w}\left(  x,t\right)   &  =\sum_{_{\mathbf{q}}\in\Omega
_{\text{generic}}}\{w_{\mathbf{q}}^{-}\mathbf{r}_{-}(\mathbf{q})\exp\left(
\lambda_{\mathbf{q}}^{-}t\right)  +w_{\mathbf{q}}^{+}\mathbf{r}_{+}%
(\mathbf{q})\exp\left(  \lambda_{\mathbf{q}}^{+}t\right)  \}e_{\mathbf{q}%
}(x)\label{l}\\
&  \ \ \ \ +\sum_{_{\mathbf{q}}\in\Omega_{\text{defective}}}\{\left(
w_{\mathbf{q}}\mathbf{r}(\mathbf{q})+w_{\mathbf{q}}^{\prime}\mathbf{r}%
^{\prime}(\mathbf{q})\right)  +w_{\mathbf{q}}^{\prime}\mathbf{r}%
(\mathbf{q})t\}\exp\left(  \lambda_{\mathbf{q}}t\right)  e_{\mathbf{q}%
}(x)\nonumber\\
&  \ \ \ \ +\sum_{_{\mathbf{q}}\in\Omega_{\text{complex}}}\{w_{\mathbf{q}%
}^{\operatorname{Re}}\left(  \operatorname{Re}\mathbf{r}(\mathbf{q}%
)\cos\left[  \left(  \operatorname{Im}\lambda_{\mathbf{q}}\right)  t\right]
-\operatorname{Im}\mathbf{r}(\mathbf{q})\sin\left[  \left(  \operatorname{Im}%
\lambda_{\mathbf{q}}\right)  t\right]  \right) \nonumber\\
&  ~~\ +w_{\mathbf{q}}^{\operatorname{Im}}\left(  \operatorname{Re}%
\mathbf{r}(\mathbf{q})\sin\left[  \left(  \operatorname{Im}\lambda
_{\mathbf{q}}\right)  t\right]  +\operatorname{Im}\mathbf{r}(\mathbf{q}%
)\cos\left[  \left(  \operatorname{Im}\lambda_{\mathbf{q}}\right)  t\right]
\right)  \}\exp[\left(  \operatorname{Re}\lambda_{\mathbf{q}}\right)
t]e_{\mathbf{q}}(x)\nonumber\\
&  \equiv e^{\mathcal{L}t}\mathbf{w}\left(  x,0\right)  .\nonumber
\end{align}
For any $\mathbf{u}\left(  \cdot\mathbf{,}t\right)  \in\left[  L^{2}\left(
{\mathbb T}^d%
\right)  \right]  ^{2}$, we denote $\left\Vert \mathbf{u}\left(
\cdot\mathbf{,}t\right)  \right\Vert \equiv\left\Vert \mathbf{u}\left(
\cdot\mathbf{,}t\right)  \right\Vert _{L^{2}}$. Our main result of this
section is

\begin{lemma}
\label{lineargrowth}Assume that (\ref{tr-det}) and the instability criterion
(\ref{instability-criteria}) are valid. Suppose
\[
\mathbf{w}\left(  x,t\right)  =[u\left(  x\mathbf{,}t\right)  ,v\left(
x\mathbf{,}t\right)  ]\equiv e^{\mathcal{L}t}\mathbf{w}\left(  x,0\right)
\]
as in (\ref{l}) is a solution to the linearized reaction-diffusion system
(\ref{linear1})-(\ref{linear2}) with initial condition $\mathbf{w}\left(
\mathbf{x},0\right)  $. Then there exists a constant $C_{1}\geq1$ depending on
$\bar{U},\bar{V},\bar{D}_{1},\bar{D}_{2},A$ such that%
\[
\left\Vert \mathbf{w}\left(  \mathbf{\cdot},t\right)  \right\Vert \leq
C_{1}\exp\left(  \lambda_{\text{max}}t\right)  \left\Vert \mathbf{w}\left(
\mathbf{\cdot},0\right)  \right\Vert ,
\]
for all $t\geq0$.
\end{lemma}

\begin{proof}
We first notice that from the quadratic formula for (\ref{dispersion}), for
$q$\ large,
\[
\left\vert \det[\mathbf{r}_{-}(\mathbf{q}),\mathbf{r}_{+}(\mathbf{q}%
)]\right\vert =\frac{\lambda_{\mathbf{q}}^{+}-\lambda_{\mathbf{q}}^{-}%
}{\left\vert \bar{f}_{v}\right\vert }\geq c\frac{\left\vert \bar{D}_{1}%
-\bar{D}_{2}\right\vert }{\left\vert \bar{f}_{v}\right\vert }q^{2}.
\]
Thus solving (\ref{inde}) yields, due to $\bar{D}_{1}\neq\bar{D}_{2},$
\begin{align*}
|w_{\mathbf{q}}^{\pm}|  &  \leq\frac{1}{\det[\mathbf{r}_{-}(\mathbf{q}%
),\mathbf{r}_{+}(\mathbf{q})]}|\mathbf{r}_{\pm}(\mathbf{q})|\times
|\mathbf{w}_{\mathbf{q}}|\\
&  \leq C|\mathbf{w}_{\mathbf{q}}|,
\end{align*}
Since $\lambda_{\mathbf{q}}<0,$ for $q\in\Omega_{\text{defective}},$ we have
\[
t\exp\left(  \lambda_{\mathbf{q}}t\right)  \leq C.
\]
Moreover, recall $\operatorname{Re}\lambda(\mathbf{q})<0$ for $q\in
\Omega_{\text{complex}}.$ Thus we deduce the Lemma on the linear growth rate
by the formula (\ref{l}).
\end{proof}

\bigskip

\section{Main Result}

Let $\theta$ be a small fixed constant, and $\lambda_{\text{max}}$ be the
dominant eigenvalue which is the maximal growth rate.\ We also denote the gap
between the largest growth rate $\lambda_{\text{max}}$ and the rest by
$\nu>0.$ Then for $\delta>0$ arbitrary small, we define the escape time
$T^{\delta}$ by%
\begin{equation}
\theta=\delta\exp\left(  \lambda_{\text{max}}T^{\delta}\right)  ,
\label{theta}%
\end{equation}
or equivalently%
\[
T^{\delta}=\frac{1}{\lambda_{\text{max}}}\ln\frac{\theta}{\delta}.
\]

Our main theorem is

\begin{theorem}
Assume (\ref{tr-det}) and that there exists $q^{2}=\sum_{i=1}^{d}q_{i}^{2}$
satisfying instability criterion (\ref{instability-criteria}). Let%
\begin{align*}
\mathbf{w}_{0}(x)  &  =\sum_{_{\mathbf{q}}\in\Omega}\{w_{\mathbf{q}}%
^{-}\mathbf{r}_{-}(\mathbf{q})+w_{\mathbf{q}}^{+}\mathbf{r}_{+}(\mathbf{q}%
)\}e_{\mathbf{q}}(x)\\
&  +\sum_{_{\mathbf{q}}\in\Omega_{\text{defective}}}\mathbf{\{}w_{\mathbf{q}%
}\mathbf{r}(\mathbf{q})+w_{\mathbf{q}}^{\prime}\mathbf{r}^{\prime}%
(\mathbf{q})\}e_{\mathbf{q}}(x)\\
&  +\sum_{_{\mathbf{q}}\in\Omega_{\text{complex}}}\mathbf{\{}w_{\mathbf{q}%
}^{\operatorname{Re}}\operatorname{Re}\mathbf{r}(\mathbf{q})+w_{\mathbf{q}%
}^{\operatorname{Im}}\operatorname{Im}\mathbf{r}(\mathbf{q})\}e_{\mathbf{q}%
}(x).
\end{align*}
$\in H^{2}$ such that $||\mathbf{w}_{0}||=1.$ Assume $D_{1},D_{2}\ f,g\in
C^{2}\mathbb{\ }$near $\bar{W},$ so that there exists $\eta>0$
\begin{equation}
C_{\eta}\equiv\max_{||w||_{\infty}\leq\eta}\{\sum_{i=1}^{2}||D_{i}(\bar
{W}+w)||_{C^{2}}+||f(\bar{W}+w)||_{C^{2}}+||g(\bar{W}+w)||_{C^{2}}<\infty.
\label{a-F}%
\end{equation}
Then there exist constants $\delta_{0}>0,$ $C>0,$ and $\theta>0,$ depending on
$\bar{U},\bar{V},\bar{D}_{1},\bar{D}_{2},f,g,$ such that for all $0<\delta
\leq\delta_{0}$, if the initial perturbation of the steady state $[\bar
{U},\bar{V}]$ in (\ref{steady-state}) is
\[
\mathbf{w}^{\delta}\left(  \mathbf{x},0\right)  =\delta\mathbf{w}_{0},
\]
then its nonlinear evolution $\mathbf{w}^{\delta}(t,x)$ satisfies
\begin{align}
&  ||\mathbf{w}^{\delta}(t,x)-\delta e^{\lambda_{\max}t}\sum_{_{\mathbf{q}}%
\in\Omega_{\text{max}}}w_{\mathbf{q}}^{+}\mathbf{r}_{+}(\mathbf{q}%
)e_{\mathbf{q}}(x)||\label{ppattern}\\
&  \leq C\{e^{-\nu t}+\delta||w_{0}||_{H^{2}}^{2}+\delta e^{\lambda_{\max}%
t}\}\delta e^{\lambda_{\max}t}\nonumber
\end{align}
for $0\leq t\leq T^{\delta},$ and $\nu>0$ is the gap between $\lambda_{\max}$
and the rest of $\operatorname{Re}\lambda_{\mathbf{q}}$ in (\ref{dispersion}).
\end{theorem}

We notice that for $0\leq t\leq T^{\delta},$ $\delta e^{\lambda_{\max}t}%
\leq\theta,\,$\ is$\ $sufficiently small. The initial profile $\mathbf{w}_{0}$
is any $H^{2}$ function. In particular, as long as $w_{\mathbf{q}_{0}}^{+}%
\neq0$ for at least one $\mathbf{q}_{0}\mathbf{\in}\Omega_{\text{max}}$
(generic for a general $H^{2}$ perturbation), the part of its fastest growing
modes satisfies
\[
||\delta e^{\lambda_{\max}t}\sum_{_{\mathbf{q}}\in\Omega_{\text{max}}%
}w_{\mathbf{q}}^{+}\mathbf{r}_{+}(\mathbf{q})e_{\mathbf{q}}||\geq\delta
e^{\lambda_{\max}t}|w_{\mathbf{q}_{0}}^{+}||\mathbf{r}_{+}(\mathbf{q}_{0})|,
\]
which has the dominant leading order of $\delta e^{\lambda_{\max}t}.$ Our
estimate (\ref{ppattern}) implies that the dynamics of a general perturbation
can be characterized by such linear dynamics over a long time period of
$\varepsilon T^{\delta}\leq t\leq T^{\delta},$ for any fixed constant
$\varepsilon>0$. In particular, choose a fixed $\mathbf{q}_{0}\in
\Omega_{\text{max}}$ and let
\[
w_{0}(x)=\frac{\mathbf{r}_{+}(\mathbf{q}_{0})}{|\mathbf{r}_{+}(\mathbf{q}%
_{0})|}e_{\mathbf{q}_{0}}(x)
\]
then if $t=T^{\delta},$
\[
\left\Vert \mathbf{w}^{\delta}(t,\cdot)-\delta e^{\lambda_{\max}T^{\delta}%
}\frac{\mathbf{r}_{+}(\mathbf{q}_{0})}{|\mathbf{r}_{+}(\mathbf{q}_{0}%
)|}e_{\mathbf{q}_{0}}(\cdot)\right\Vert \leq C\{\delta^{\nu/\lambda_{\text{max
}}}+\theta^{2}\},
\]
hence
\[
\left\Vert \mathbf{w}^{\delta}(t,\cdot)\right\Vert \geq\theta-C\{\delta
^{\nu/\lambda_{\text{max }}}+\theta^{2}\}\geq\theta/2>0,
\]
which implies nonlinear instability as $\delta\rightarrow0$. The instability
occurs before the possible blow-up time.

Reaction-diffusion systems are often employed to study chemical and biological
pattern formation and have received much attention from scientists \cite{DMO},
\cite{Gr}, \cite{M}, \cite{Me}, \cite{Ty}, since the pioneering work of Turing
\cite{T} in 1951. This symmetry breaking instability is called
diffusion-driven instability, since the presence of diffusion and the
difference of diffusion coefficients are essential for the instability
mechanism and nonuniform pattern formation. After some experimental results
such as in \cite{CDBK}, \cite{LMPS}, \cite{QS}, more extensive and serious
works began towards this Turing-like pattern formation across many fields of
study. Our result can be interpreted as a mathematical description of early
pattern formation. Each initial perturbation can be drastically different from
another, which gives rise to the richness of the pattern; on the other hand,
the finite number maximal growing modes determine the common characteristics
of the pattern, over the time scale of $\ln\frac{1}{\delta}.$ In comparision
with an earlier \textit{different} result along this direction \cite{ST}:
First of all, the reaction-diffusion system considered here is \textit{not}
scaled. Secondly, our initial perturbation is more general, need not be close
to the space of finite number of maximal growing modes. Thirdly, a precise
estimate of the time scale ($\ln\frac{1}{\delta}$) for pattern formation is
given here, without an \textit{a-priori} assumption for the smallness of the
perturbation later in time as in \cite{ST}. Lastly, based on Guo-Strauss'
bootstrap argument, our proof is much simpler and direct.

\section{ Bootstrap Lemma}

We state existence of local-in-time solutions for (\ref{nonlinear1}%
)-(\ref{nonlinear2}).

\begin{lemma}
(Local existence) For $s\geq1$ $\left(  d=1\right)  $ and $s\geq2$ $\left(
d=2,3\right)  $, there exist a $T>0$ and a constant $C$ depending on $\bar
{U},\bar{V},D_{1},D_{2},f,g$ such that $\left\Vert \mathbf{w}(t)\right\Vert
_{H^{s}}$ is continuous in $[0,T),$ and
\[
\left\Vert \mathbf{w}(t)\right\Vert _{H^{s}}\leq C\left\Vert \mathbf{w}\left(
0\right)  \right\Vert _{H^{s}}.
\]

\end{lemma}

We now derive the following energy estimates for $d$-dimensional
reaction-diffusion system with $d=1,2,3$.

\begin{lemma}
Suppose that $[u\left(  x\mathbf{,}t\right)  ,v\left(  x,t\right)  ]$ is a
solution to the full system (\ref{nonlinear1})-(\ref{nonlinear2}). Then for
$||w(t)||_{H^{2}}\leq\eta,$
\begin{align*}
&  \frac{1}{2}\frac{d}{dt}\sum_{\left\vert \partial\right\vert =2}%
\int_{\mathbb{T}^{d}}\{\left\vert \partial u\right\vert ^{2}+\left\vert
\partial v\right\vert ^{2}\}d\mathbf{x}\\
&  +\sum_{\left\vert \partial\right\vert =2}\int_{\mathbb{T}^{d}}\left\{
\frac{\bar{D}_{1}}{2}\left\vert \nabla\partial u\right\vert ^{2}+\bar{D}%
_{2}\left\vert \nabla\partial v\right\vert ^{2}\right\}  d\mathbf{x}%
+\frac{\left\vert \bar{g}_{v}\right\vert }{2}\sum_{\left\vert \alpha
\right\vert =2}\int_{\mathbb{T}^{d}}|\partial v|^{2}\\
&  \leq C_{0}C_{1}||\mathbf{w}||_{H^{2}}||\nabla^{3}\mathbf{w}||^{2}%
+C_{2}||u\mathbf{||}^{2}.
\end{align*}
where $C_{0}$ is the universal constant while $C_{1}=C_{0}C_{\eta}(1+\eta)$
and
\[
C_{2}=\frac{\left(  \frac{\left(  \bar{f}_{v}+\bar{g}_{u}\right)  ^{2}%
}{2\left\vert \bar{g}_{v}\right\vert }+\bar{f}_{u}\right)  ^{3}}{\bar{D}%
_{1}^{2}}.
\]

\end{lemma}

\begin{proof}
We first notice that the reaction-diffusion system (\ref{nonlinear1}%
)-(\ref{nonlinear2}) preserves the evenness of the solution $\mathbf{w}(x,t),$
i.e., if $\mathbf{w}(x,t)$ is a solution, then $\mathbf{w}(-x_{i},t)$ is also
a solution. We can regard the Neumann problem as a special case with evenness
of the periodic problem by standard way of even extension $\mathbf{w}(x,t)$
with respect to one of the $x_{i}$. For this reason we may assume periodicity
at the boundary of the extended periodic box $2\mathbf{T}^{3}\equiv(-\pi
,\pi)^{d}$. Since now there is no contributions from the boundaries, we can
take second order $\partial$-derivative of (\ref{matrix}) to get%
\begin{equation}
\frac{1}{2}\frac{d}{dt}\int_{2\mathbb{T}^{d}}\left\vert \partial
\mathbf{w}\right\vert ^{2}=\int_{2\mathbb{T}^{d}}\partial\mathbf{w}%
^{T}\partial\mathcal{L}\left(  \mathbf{w}\right)  \mathbf{+}\int
_{2\mathbb{T}^{d}}\partial\mathbf{w}^{T}\partial\mathcal{N}\left(
\mathbf{w}\right)  . \label{high}%
\end{equation}

We first treat the last nonlinear term:
\begin{align*}
&  -\int_{2\mathbb{T}^{d}}\{\nabla\partial\mathbf{w\}}^{T}[\partial\{D\left(
\mathbf{w}+\mathbf{\bar{W}}\right)  \nabla\mathbf{w}\}+\bar{D}\nabla
\partial\mathbf{w}]+\{\nabla\partial\mathbf{w\}}^{T}\partial\left(
\mathbf{F-}A\mathbf{w}\right) \\
\leq &  C\left\Vert D\left(  \mathbf{w}+\mathbf{\bar{W}}\right)  -\bar
{D}\right\Vert _{\infty}\left\Vert \nabla\partial\mathbf{w}\right\Vert
^{2}+C\left\Vert \left(  \nabla D\right)  \left(  \mathbf{w}+\mathbf{\bar{W}%
}\right)  \right\Vert _{\infty}\left\Vert \nabla\mathbf{w}\right\Vert
_{\infty}\left\Vert \partial\mathbf{w}\right\Vert \left\Vert \nabla
\partial\mathbf{w}\right\Vert \\
&  +C\left\Vert \left(  \partial D\right)  \left(  \mathbf{w}+\mathbf{\bar{W}%
}\right)  \right\Vert _{\infty}\left\Vert \nabla\mathbf{w}\right\Vert _{L^{4}%
}^{2}\left\Vert \nabla\mathbf{w}\right\Vert _{\infty}\left\Vert \nabla
\partial\mathbf{w}\right\Vert \\
&  +C\left\Vert \left(  \partial\mathbf{F}\right)  \left(  \mathbf{w}%
+\mathbf{\bar{W}}\right)  \right\Vert _{\infty}\left\Vert \nabla
\mathbf{w}\right\Vert _{\infty}\left\Vert \nabla\mathbf{w}\right\Vert
\left\Vert \nabla\partial\mathbf{w}\right\Vert +C\left\Vert \nabla
\mathbf{F}\left(  \mathbf{w}+\mathbf{\bar{W}}\right)  \mathbf{-}A\right\Vert
_{\infty}\left\Vert \partial\mathbf{w}\right\Vert \left\Vert \nabla
\partial\mathbf{w}\right\Vert .
\end{align*}
We apply the following the Sobolev imbedding to control $||w||_{\infty}$
\begin{equation}
\left\Vert g\right\Vert _{L^{\infty}\left(  2\mathbb{T}^{d}\right)  }\leq
C_{0}\left\Vert g\right\Vert _{H^{2}\left(  2\mathbb{T}^{d}\right)  },
\label{Sobolev1}%
\end{equation}
$\ $\ for $d\leq3.$ Moreover, from the periodic boundary conditions,
\[
\int_{2\mathbb{T}^{d}}\nabla u=\int_{2\mathbb{T}^{d}}\nabla v=0,
\]
we also use the Poincare inequality%
\begin{equation}
||g||\leq\left\Vert g\right\Vert _{L^{4}(2\mathbb{T}^{d})}\leq C_{0}\left\Vert
\nabla g\right\Vert \text{ \ \ \ if }d\leq3. \label{Sobolev2}%
\end{equation}
to further get%
\[
||\nabla\mathbf{w}||_{\infty}\leq C_{0}\left\Vert \nabla\mathbf{w}\right\Vert
_{H^{2}}\leq C_{0}\sum_{|\partial|=2}||\partial\nabla\mathbf{w}||.
\]
where $C_{0}$ is a universal constant. From (\ref{a-F}) and the assumption
$||\mathbf{w}||_{H^{2}}\leq\eta$, the last nonlinear term in (\ref{high}) is
bounded by
\[
C_{0}C_{\eta}(1+\eta)\left\Vert \mathbf{w}\right\Vert _{H^{2}}\left\Vert
\nabla\partial\mathbf{w}\right\Vert ^{2}.
\]
We now estimate the second quadratic term in (\ref{high})
\begin{align*}
&  -\int_{2\mathbb{T}^{d}}\{\bar{D}_{1}\left\vert \nabla\partial u\right\vert
^{2}+\bar{D}_{2}\left\vert \nabla\partial v\right\vert ^{2}\}+\bar{g}_{v}%
\int_{2\mathbb{T}^{d}}|\partial v|^{2}\\
&  +\left(  \bar{f}_{v}+\bar{g}_{u}\right)  \int_{2\mathbb{T}^{d}}\partial
u\partial v+\bar{f}_{u}\int_{2\mathbb{T}^{d}}\left\vert \partial u\right\vert
^{2}.
\end{align*}

The last two terms are bounded by%
\begin{align*}
&  \left(  \bar{f}_{v}+\bar{g}_{u}\right)  \int_{2\mathbb{T}^{d}}\partial
u\partial v+\bar{f}_{u}\int_{2\mathbb{T}^{d}}\left\vert \partial u\right\vert
^{2}\\
&  \leq\frac{\left\vert \bar{g}_{v}\right\vert }{2}\int_{2\mathbb{T}^{d}%
}\left\vert \partial v\right\vert ^{2}+\{\frac{\left(  \bar{f}_{v}+\bar{g}%
_{u}\right)  ^{2}}{2\left\vert \bar{g}_{v}\right\vert }+\bar{f}_{u}%
\}\int_{2\mathbb{T}^{d}}\left\vert \partial u\right\vert ^{2}.
\end{align*}

Thus we can bound the linear term in (\ref{high}) by $(\bar{g}_{v}<0)$%
\begin{align*}
&  -\int_{2\mathbb{T}^{d}}\{\bar{D}_{1}\left\vert \nabla\partial u\right\vert
^{2}+\bar{D}_{2}\left\vert \nabla\partial v\right\vert ^{2}\}-\frac{\left\vert
\bar{g}_{v}\right\vert }{2}\int_{2\mathbb{T}^{d}}\left\vert \partial
v\right\vert ^{2}\\
&  +\{\frac{\left(  \bar{f}_{v}+\bar{g}_{u}\right)  ^{2}}{2\left\vert \bar
{g}_{v}\right\vert }+\bar{f}_{u}\}\int_{2\mathbb{T}^{d}}\left\vert \partial
u\right\vert ^{2}.
\end{align*}

$\ \ \ \ \ \ $By the interpolation between $\left\Vert \nabla\partial
u\right\Vert $ and $||u\mathbf{||}$ , the last term above is bounded by
\[
\{\frac{\left(  \bar{f}_{v}+\bar{g}_{u}\right)  ^{2}}{2\left\vert \bar{g}%
_{v}\right\vert }+\bar{f}_{u}\}\{a\int_{2\mathbb{T}^{d}}\left\Vert
\nabla\partial u\right\Vert ^{2}+\frac{1}{4a^{2}}\int_{2\mathbb{T}^{d}%
}||u\mathbf{||}^{2}\}
\]
for any $a>0.$ We can choose $a$ such that
\[
\{\frac{\left(  \bar{f}_{v}+\bar{g}_{u}\right)  ^{2}}{2\left\vert \bar{g}%
_{v}\right\vert }+\bar{f}_{u}\}a=\frac{1}{2}\bar{D}_{1}.
\]
Collecting terms, we conclude the proof.
\end{proof}

We are now ready to establish the bootstrap lemma, which controls the $H^{2}$
growth of $\mathbf{w}(x,t)$ in term of its $L^{2}$ growth nonlinearly.

\begin{lemma}
\label{wh2}Suppose that $\mathbf{w}(x,t)$ is a solution to the full system
(\ref{nonlinear1})-(\ref{nonlinear2}) such that for $0\leq t\leq T$
\[
||\mathbf{w}(\cdot,t)||_{H^{2}}\leq\min\left\{  \eta,\frac{\bar{D}_{1}}%
{2C_{0}C_{1}},\frac{\bar{D}_{2}}{C_{0}C_{1}}\right\}
\]
and
\begin{equation}
||\mathbf{w}(\cdot,t)||\leq2C_{1}e^{\lambda_{\max}t}||\mathbf{w}(\cdot,0)||,
\label{l2growth}%
\end{equation}
then we have for $0\leq t\leq T$
\[
||\mathbf{w}(t)||_{H^{2}}^{2}\leq C_{3}\{||\mathbf{w}(0)||_{H^{2}}%
^{2}+e^{2\lambda_{\max}t}||\mathbf{w}(\cdot,0)||^{2}\}
\]
where $C_{3}$ $=C_{1}^{2}\max\{\frac{4C_{2}}{\lambda_{\text{max}}},1\}\geq1.$
\end{lemma}

\begin{proof}
It suffices to only consider the second-order derivatives of $\mathbf{w}%
(x,t).$ From the previous lemma and our assumption for $||\mathbf{w}||_{H^{2}%
},$ we deduce that
\[
\frac{1}{2}\frac{d}{dt}\sum_{\left\vert \alpha\right\vert =2}\int
_{\mathbb{T}^{d}}\left\{  |\partial u|^{2}+|\partial v|^{2}\right\}
d\mathbf{x}\leq C_{2}||u\mathbf{||}^{2}.
\]
So that by (\ref{l2growth}) and an integration from $0$ to $t\leq T,$ we have
\begin{align*}
&  \sum_{|\partial|=2}\int_{\mathbb{T}^{d}}\left\{  |\partial u(t)|^{2}%
+|\partial v(t)|^{2}\right\} \\
&  \leq\sum_{|\partial|=2}\int_{\mathbb{T}^{d}}\left\{  |\partial
u(0)|^{2}+|\partial v(0)|^{2}\right\}  +\frac{4C_{2}C_{1}^{2}}{\lambda
_{\text{max}}}e^{2\lambda_{\max}t}||\mathbf{w}(\cdot,0)||^{2}.
\end{align*}
Thus our lemma follows.
\end{proof}

\section{Nonlinear instability and pattern formation}

We now prove our main Theorem 1:

\begin{proof}
Let $\mathbf{w}^{\delta}\left(  x,t\right)  $ be the family of solutions to
the reaction-diffusion system (\ref{nonlinear1})-(\ref{nonlinear2}) with
initial data $\mathbf{w}^{\delta}\left(  x,0\right)  =\delta\mathbf{w}_{0}$.
Define $T^{\ast}$ by%
\[
T^{\ast}=\sup\left\{  t\ |~\left\Vert \mathbf{w}^{\delta}(t)-\delta
e^{\mathcal{L}t}\mathbf{w}_{0}\right\Vert \leq\frac{C_{1}}{2}\delta\exp\left(
\lambda_{\max}t\right)  \right\}  .
\]
Note that $T^{\ast}$ is well defined. We also define%
\[
T^{\ast\ast}=\sup\left\{  t\ |~||\mathbf{w}(t)||_{H^{2}}\leq\min\left\{
\eta,\frac{\bar{D}_{1}}{2C_{0}C_{1}},\frac{\bar{D}_{2}}{C_{0}C_{1}}\right\}
\right\}  .
\]
\ 

\ \ \ \ \ We now derive estimates for $H^{2}$ norm of $\mathbf{w}^{\delta
}(x,t)$ for $0\leq t\leq\min\{T^{\ast},T^{\ast\ast}\}.$ First of all, by the
definition of $T^{\ast},$ for $t\leq T^{\ast}$ and Lemma \ref{lineargrowth}
\[
\left\Vert \mathbf{w}^{\delta}(t)\right\Vert \leq\frac{3C_{1}}{2}\delta
\exp\left(  \lambda_{\max}t\right)  .
\]
\ \ \ \ \ \ \ 

Moreover, using Lemma 4 and applying a bootstrap argument yields
\begin{equation}
\left\Vert \mathbf{w}^{\delta}(t)\right\Vert _{H^{2}}\leq\sqrt{C_{3}}%
\{\delta||\mathbf{w}_{0}||_{H^{2}}+\delta e^{\lambda_{\max}t}\}. \label{h2}%
\end{equation}
\ 

We now estimate the $L^{2}$ norm of $\mathbf{w}^{\delta}(x,t)$ for $0\leq
t\leq\min\{T^{\ast},T^{\ast\ast}\}.$ We apply Duhamel's principle to obtain%
\[
\mathbf{w}^{\delta}\left(  t\right)  =\delta e^{\mathcal{L}t}\mathbf{w}%
_{0}-\int_{0}^{t}e^{\mathcal{L}\left(  t-\tau\right)  }\mathcal{N}\left(
\mathbf{w}^{\delta}\left(  \tau\right)  \right)  d\tau,
\]
Using Lemma 1, (\ref{Sobolev1}), (\ref{Sobolev2}), and Lemma 4 yields, for
$0\leq t\leq\min\{T^{\ast},T^{\ast\ast}\}$%
\begin{align*}
&  \left\Vert \mathbf{w}^{\delta}\left(  t\right)  -\delta e^{\mathcal{L}%
t}\mathbf{w}_{0}\right\Vert \\
\leq &  C_{1}\int_{0}^{t}e^{\lambda_{\max}\left(  t-\tau\right)  }\left\Vert
\mathbf{\{}\nabla\cdot\left(  D\nabla\mathbf{w}^{\delta}\right)  -\bar
{D}\nabla^{2}\mathbf{w^{\delta}\}+F-}A\mathbf{w}^{\delta}\right\Vert d\tau\\
\leq &  C_{1}\int_{0}^{t}e^{\lambda_{\max}\left(  t-\tau\right)  }%
||D||_{C^{1}}\left\Vert \mathbf{w}^{\delta}\left(  \tau\right)  \right\Vert
_{\infty}\left\Vert \mathbf{w}^{\delta}\left(  \tau\right)  \right\Vert
_{H^{2}}d\tau\\
&  +C_{1}\int_{0}^{t}e^{\lambda_{\max}\left(  t-\tau\right)  }||D||_{C^{1}%
}\left\Vert \nabla\mathbf{w}^{\delta}\left(  \tau\right)  \right\Vert _{L^{4}%
}\left\Vert \nabla\mathbf{w}^{\delta}\left(  \tau\right)  \right\Vert _{L^{4}%
}d\tau\\
&  +C_{1}\int_{0}^{t}e^{\lambda_{\max}\left(  t-\tau\right)  }||F||_{C^{2}%
}\left\Vert \mathbf{w}^{\delta}\left(  \tau\right)  \right\Vert _{\infty
}\left\Vert \mathbf{w}^{\delta}\left(  \tau\right)  \right\Vert d\tau\\
\leq &  C_{1}C_{0}^{2}C_{\eta}\int_{0}^{t}e^{\lambda_{\max}\left(
t-\tau\right)  }\left\Vert \mathbf{w}^{\delta}\left(  \tau\right)  \right\Vert
_{H^{2}}^{2}d\tau.
\end{align*}
from assumption (\ref{a-F}) with $||w||_{H^{2}}\leq\eta.$ We plug (\ref{h2})
with $t=\tau$ to further obtain
\begin{align}
&  \left\Vert \mathbf{w}^{\delta}\left(  t\right)  -\delta e^{\mathcal{L}%
t}\mathbf{w}_{0}\right\Vert \label{l2}\\
&  \leq C_{1}C_{0}^{2}C_{\eta}C_{3}\int_{0}^{t}e^{\lambda_{\max}\left(
t-\tau\right)  }\{\delta^{2}||\mathbf{w}_{0}||_{H^{2}}^{2}+\delta
^{2}e^{2\lambda_{\max}\tau}\}d\tau\nonumber\\
&  \leq C_{1}C_{0}^{2}C_{\eta}C_{3}\{\frac{||\mathbf{w}_{0}||_{H^{2}}%
^{2}\delta}{\lambda_{\text{max}}}+\frac{1}{\lambda_{\text{max}}}\delta
e^{\lambda_{\max}t}\}\delta e^{\lambda_{\max}t}.\nonumber
\end{align}
\ \ \ \ \ \ 

We now choose $\theta$ in $T^{\delta}$ in (\ref{theta}) to satisfy
\begin{align}
C_{0}^{2}C_{3}C_{\eta}\theta &  <\frac{\lambda_{\text{max}}}{4},\label{th1}\\
2\sqrt{C_{3}}\theta &  <\min\left\{  \eta,\frac{\bar{D}_{1}}{2C_{0}C_{1}%
},\frac{\bar{D}_{2}}{C_{0}C_{1}}\right\}  . \label{th2}%
\end{align}
We now prove by contradiction that for $\delta$ sufficiently small,
\[
T^{\delta}\leq\min\{T^{\ast},T^{\ast\ast}\},
\]
and therefore our theorem follows from (\ref{l2}), by further separating
$\mathbf{q}\in\Omega_{\max}$ and move $\mathbf{q}\notin\Omega_{\max}$ in
(\ref{l}) to the right hand side .

\ \ \ If $T^{\ast\ast}$ is the smallest among $T^{\delta},$ $T^{\ast}$ and
$T^{\ast\ast}$, we can let $t=T^{\ast\ast}<T^{\delta}$ in (\ref{h2})
\begin{align*}
\left\Vert \mathbf{w}^{\delta}(T^{\ast\ast})\right\Vert _{H^{2}}  &
<\sqrt{C_{3}}\{\delta||\mathbf{w}_{0}||_{H^{2}}+\delta e^{\lambda_{\max
}T^{\delta}}\}\\
&  =\sqrt{C_{3}}\{\delta||\mathbf{w}_{0}||_{H^{2}}+\theta\}\leq2\sqrt{C_{3}%
}\theta,
\end{align*}
for small $\delta$ such that $\delta||\mathbf{w}_{0}||_{H^{2}}\leq\theta$. By
the choice of $\theta$ in (\ref{th2}), we have%
\[
||\mathbf{w}(T^{\ast\ast})||_{H^{2}}<\min\left\{  \eta,\frac{\bar{D}_{1}%
}{2C_{0}C_{1}},\frac{\bar{D}_{2}}{C_{0}C_{1}}\right\}  .
\]
This is a contradiction to the definition of $T^{\ast\ast}.$

\ \ \ On the other hand, if $T^{\ast}$ is the smallest among among $T^{\delta
},$ $T^{\ast}$ and $T^{\ast\ast}$, we can let $t=T^{\ast}$ in (\ref{l2}) to
get%
\begin{align*}
&  \left\Vert \mathbf{w}^{\delta}\left(  T^{\ast}\right)  -\delta
e^{\mathcal{L}t}\mathbf{w}_{0}\right\Vert \\
&  \leq C_{1}C_{0}^{2}C_{3}C_{\eta}\{\frac{||\mathbf{w}_{0}||_{H^{2}}%
^{2}\delta}{\lambda_{\text{max}}}+\frac{1}{\lambda_{\text{max}}}\delta
e^{\lambda_{\max}T^{\delta}}\}\delta e^{\lambda_{\max}T^{\ast}}\\
&  \leq C_{1}C_{0}^{2}C_{3}C_{\eta}\{\frac{||\mathbf{w}_{0}||_{H^{2}}%
^{2}\delta}{\lambda_{\text{max}}}+\frac{\theta}{\lambda_{\text{max}}}\}\delta
e^{\lambda_{\max}T^{\ast}}\\
&  <\frac{C_{1}}{2}\delta e^{\lambda_{\max}T^{\ast}},
\end{align*}
for $C_{0}^{2}C_{3}C_{\eta}\frac{||\mathbf{w}_{0}||_{H^{2}}^{2}\delta}%
{\lambda_{\text{max}}}<1/4$ for $\delta$ small, by our choice of $\theta$ in
(\ref{th1}). This again contradicts the definition of $T$ and our theorem follows.
\end{proof}

\end{document}